 \documentstyle[12pt]{article}

\newcounter{const}

\newcommand{\argm}{\mathop{\rm argmax}\limits}

\begin{document}




\begin{center}
\vspace{2mm}
{\sc SCHL\"OMILCH AND BELL SERIES FOR BESSEL'S
 FUNCTIONS, WITH PROBABILISTIC APPLICATIONS. }\\

\vspace{4mm}

{\sc Ostrovsky E. and Sirota L.} \footnote{Corresponding author. Department of Mathematics,
 Bar - Ilan University, Ramat Gan, Israel, 59200.\\
 E - mail: sirota@zahav.net.il}\\

Department of Mathematics, Bar - Ilan University, Ramat - Gan, 59200,Israel.\\
e - mail: \ sirota@zahav.net.il\\

\end{center}

\vspace{3mm}

 {\sc \textbf{Abstracts}.} \par
 \vspace{2mm}
 We have introduced and investigated so-called Schl\"omilch's and
Bell's series for modified Bessel's functions, namely, their
asymptotic and non-asymptotic properties, connection with
Stirling's and Bell's numbers etc.\par
  We have obtained exact constants in the moment
inequalities for sums of centered independent random variables,
improved their asymptotical properties, found lower and upper
bounds, calculated a more exact approximation, elaborated the
numerical algorithm for their calculation, studied the class of
smoothing, etc.


\vspace{3mm}


\vspace{2mm}


\vspace{2mm}

\section{Introduction. Statement of problem.} \par
 Let $ I_{\nu}(z) $ be usually {\bf modified } Bessel's function:

$$
I_{\nu}(z) = 2^{- \nu} \ z^{\nu} \ \sum_{k=0}^{\infty}
4^{-k} z^{2k} /(k! \ \Gamma(\nu+k+1) ), \eqno(1.0)
$$
 where $ \Gamma(\cdot) $ denotes the gamma-function and $ \nu \ge 0. $ \par
 We define the so-called Schl\"omilch's functions $ F_3, F_2, F_1 $
{\it of the first kind} for the values $ p \ge 0, \ \theta > 0, \
\beta > 0: $

$$
F_3^{(S)} = F_3(p;\theta,\beta) =
 \sum_{k = - \infty}^{\infty} |k|^p \ \theta^k \ I_k(\beta), \eqno(1.1a)
$$

 $$
   F_2(p;\beta) = F_3(p;1,\beta) = \sum_{k=-\infty}^{\infty} |k|^p \
  I_k(\beta), \eqno(1.1b)
 $$

 $$
  F_1(p) = F_3(p,1,1) = \sum_{k= - \infty}^{\infty} |k|^p \ I_k(1),
 \eqno(1.1c)
  $$

 We also define  the Schl\"omilch's functions {\it of a second kind}
$ G_3, G_2, G_1 $  for positive {\it integer} values $ p = 1,2,3,
\ldots $ only and $ \theta > 0, \ \beta > 0 $ as follows:

$$
G_3^{(S)} = G_3(p;\theta,\beta) =
 \sum_{k = - \infty}^{\infty} k^p \ \theta^k \ I_k(\beta), \eqno(1.1d)
$$

 $$
 G_2(p;\beta) = G_3(p,1,\beta) = \sum_{k = -\infty}^{\infty} k^p \ I_k(\beta),
\eqno(1.1e)
 $$

  $$
  G_1(p)) = G_3(p,1,1) = \sum_{k = -\infty}^{\infty} k^p \ I_k(1). \eqno(1.1f)
   $$
 Of course, $ G_3(2m; \theta,\beta) = F_3(2m; \theta, \beta),
 \ m = 1,2,3, \ldots.$ \par

 Recall that $ I_{-k}(z) = I_k(z), \ k = 1,2,3,... $  Therefore,

$$
F_3(p; \theta, \beta) = F_3(p; 1/\theta, \beta).
$$

  As for the definition in the case {\it usual, i.e. non-modified} Bessel's
  function and respective preliminary results see, for example,  [1], [2], p.
  38 - 42. \par

 Let us examine the probabilistic meaning of these functions. Let
$ \xi, \ \eta $ be independent random variables characterized by
the Poisson distribution  with parameters $ \lambda > 0 $ and $
\mu
> 0 $ correspondently: $ Law(\xi) = Poisson(\lambda), \ Law(\eta)
= Poisson(\mu), $  or, in more detail:

$$
{\bf P}(\xi = k) = \exp(-\lambda) \ \lambda^k/k!, \ k = 0,1,2, \ldots;
$$

$$
{\bf P}(\eta = l) = \exp(-\mu) \ \mu^l/l!, \ l = 0,1,2, \ldots;
$$
and  denote $ \tau = \xi - \eta. $  Then we obtain for
non-negative integer values $ n: $

$$
{\bf P}(\tau = n) = \sum_{k=0}^{\infty} {\bf P}(\xi = n+k) \ {\bf P}
(\eta = k) =
$$

$$
\sum_{k = 0}^{\infty} \exp(- \lambda - \mu) \ \lambda^{n + k} \ \mu^k \
/( (k+n)! \ n!) =
$$

$$
\exp(-(\lambda + \mu)) \ (\lambda/\mu))^{n/2} \
 I_n \left(2 \ \sqrt{\lambda \ \mu} \right), \ n = 0, 1, 2, \ldots
$$
 and therefore
 $$
 {\bf E}|\tau|^p = \exp(-(\lambda + \mu)) \sum_{n=- \infty}^{\infty}
 |n|^p \ (\lambda/\mu)^{n/2} \ I_n \left( 2 \ \sqrt{\lambda \ \mu } \right) =
$$

$$
\exp(-(\lambda + \mu)) \ F_3 \left(p,\sqrt{\lambda/\mu}, 2 \ \sqrt{\lambda \ \mu } \right), \eqno(1.2)
$$

 $$
 {\bf E}\tau^p = \exp(-(\lambda + \mu))\sum_{n=- \infty}^{\infty} n^p \
 (\lambda/\mu)^{n/2} \ I_n( 2 \ \sqrt{\lambda \ \mu)} =
$$

$$
\exp(-(\lambda + \mu)) \ G_3 \left(p,\sqrt{\lambda/\mu}, 2 \
\sqrt{\lambda \ \mu } \right). \eqno(1.3)
$$

  (The case of negative integers  $ n $ is considered in a similar way.) \par
 If $ \lambda = \mu, $ then
$$
{\bf E}|\tau|^p = \exp(-2 \ \lambda) \ F_2(p,2 \ \lambda),
$$
and in the case of  $ \lambda = \mu = 1/2: $ we obtain:
$$
{\bf E}|\tau|^p = F_1(p)/e.
$$

 Here we introduce some new functions which represent generalizations
 of the classical Bell's number and functions.  \par

Generalized Bell functions {\it of a first kind} $ B_4(p;a,
\lambda,\gamma)) $ can be defined as

$$
B_4(p; a,\lambda, \gamma) = \sum_{n= 0}^{\infty} \frac{|n-a|^p
\ \lambda^n }{\exp(\lambda) \ \cdot \Gamma(n + \gamma + 1)}, \eqno(1.3a)
$$

where $ \lambda > 0, \gamma \ge 0. $ We also define

 $$
 B_3(p;a,\lambda) = B_4(p;a,\lambda,0) = \sum_{n=0}^{\infty}
\frac{|n-a|^p \lambda^n}{\exp(\lambda) \ n!}, \eqno(1.3b)
 $$

 $$
 B_2(p; \lambda) = B_3(p; 0,\lambda); \eqno(1.3c)
 $$

 $$
 B(p) = B_1(p) =
  B_2(p;1) = \sum_{n=0}^{\infty} \frac{n^p}{e \ n!}. \eqno(1.3d)
$$

 We also define generalized Bell functions {\it of a second  kind}
$ D_4(p;a, \lambda,\gamma)) $  for integer positive   $ p $ values
as follows:

$$
D_4(p; a,\lambda, \gamma) = \sum_{n= 0}^{\infty} \frac{(n-a)^p \ \lambda^n }
{\exp(\lambda) \ \cdot \Gamma(n + \gamma + 1)}, \eqno(1.4a)
$$
$ \lambda > 0, \gamma \ge 0. $ We also define  $ D_3(p;a,\lambda)
= D_4(p;a,\lambda,0); \ D_2(p; \lambda) = D_3(p; 0,\lambda); \
D_1(p) = D_2(p;1) = B_1(p). $

 The probabilistic interpretation is
 $$
{\bf E}|\xi - a|^p = e^{-\lambda} \ \sum_{n=0}^{\infty}|n-a|^p
\lambda^n/n! = B_3(p;a,\lambda)
$$
 and $ \ {\bf E}(\xi - a)^p = D_3(p;a,\lambda). $ \par

 {\bf In this paper we investigate the estimations and asymptotic for
introduced functions as $ p \to \infty. $ \par
 In Section 4 we find exact constants for
moment estimations of sums of independent random variables.}\par

 \section{Non-asymptotic properties. Relations with Stirling's numbers.}\par

  Recall that the Stirling numbers of the second kind $ \{ s(n,r) \} $
 which have appeared in  the combinatorial theory ([3], p. 117, [4],
 p. 234 - 243), are defined by the following identities:

 $$
 x^n = \sum_{r=0}^n s(n,r) x_{(r)};  x_{(r)} = x(x-1)(x-2)\ldots (x-r+1),
  \ x_{(0)}=1.
 $$

 {\bf 1}. Define the function
 $ Q_{2m}(\lambda,\mu) $ for the values $ m = 1,2,3, \ldots, \
\lambda,\mu > 0  $ by the following way:

$$
Q_{2m}(\lambda, \mu) \stackrel{def}{=} \exp( - \lambda - \mu) \
F_3 \left(2m; \sqrt{\lambda/\mu}, 2 \sqrt{\lambda \mu} \right).
$$
 Proposition:
 $$
 Q_{2m}(\lambda, \ \mu) = \sum_{r=0}^{2m} \sum_{i=0}^r \sum_{j=0}^{2m - r} (-1)^r
{2m \choose r} \ \lambda^i \ \mu^j  \ s(r,i) \ s(2m - r,j).
 $$
 Therefore, the function $ Q_{2m}(\lambda, \mu) $ is a polynomial of the power
$ 2m. $ \par
 {\bf Proof}. As known  $ Q_{2m}(\lambda,\mu) = {\bf E}(\xi - \eta)^{2m}.  $
 The conclusion {\bf 1}  follows from the definition of Stirling's
 numbers, binomial formula and the following identity

$$
{\bf E} \xi_{(r)} = \lambda^r.
$$
 {\bf 2.} Let us consider the function $ D_3(p; a,\lambda) $ for integer
non-negative values $ p. $ We have:

 $$
 D_3(p; a,\lambda) = {\bf E}(\xi - a)^p = \sum_{l=0}^p (-1)^l {p \choose l} \
 a^{p-l} \sum_{r=0}^l \lambda^r s(l,r) =
 $$

 $$
\sum_{r=0}^p \lambda^r  \sum_{l=r}^p (-1)^l {p \choose l} a^{p-l} s(l,r).
 $$

 {\bf 3.} Let us denote for integer positive $ p $ the values:
  $$
   H(p; \lambda,\mu) = \exp( - \lambda - \mu) \ G_3(p; \sqrt{\lambda/\mu},
  2 \ \sqrt{\lambda \mu} ),
  $$
 where $ \lambda = \lambda_1 + \lambda_2, \ \mu = \mu_1 + \mu_2, \ \lambda_{1,2} > 0,
\ \mu_{1,2} > 0. $ We assert that

 $$
 H(p; \lambda_1 + \lambda_2, \mu_1 + \mu_2) =
 \sum_{r=0}^p  {p \choose r} H(r; \lambda_1, \mu_1) \ H(p-r; \lambda_2, \mu_2).
 $$

  Namely, we can represent (on some Probability space)
 the random variables $ \xi $ and $ \eta $ as a sums
 $ \xi = \xi_1 + \xi_2, \ \eta = \eta_1 + \eta_2, $ where all the variables
 are independent, $ Law(\xi_i) = Poisson(\lambda_i), \ Law(\eta_i) =
 Poisson(\mu_i), \ i = 1,2; $ then

 $$
   H(p; \lambda_1 + \lambda_2, \mu_1 + \mu_2) = {\bf E} \left( (\xi_1 -\eta_1) +
 (\xi_2 - \eta_2) \right)^p =
  \sum_{r=0}^p {p \choose r} {\bf E}(\xi_1 - \eta_1)^r \
 {\bf E}(\xi_2 - \eta_2)^{p - r} =
  $$

  $$
  \sum_{r=0}^p {p \choose r} H(r; \lambda_1, \mu_1) \ H(p-r, \lambda_2, \mu_2).
  $$

 \section{Asymptotic results.}\par

  Denote $ q = q(p; a,\gamma,\lambda) = (p - a - \gamma - 0.5)/\lambda, $
 $$
 p \ge 4 \ \Rightarrow g(p) = \frac{p}{e \log p}, \Delta(p) =
  \frac{\log \log p}{ \log p},
 $$
 and assume $ a,\gamma,\lambda = const, \ p \to \infty. $ \par
  Let us introduce the following relations of equivalence:
 for two positive functions $ f_1(p), f_2(p) $ we will write
 $ \ f_1 \sim f_2 $ iff

 $$
  f_1(p) = f_2(p)(1 + O(1/\log p)).
 $$
  We can also write  $ f_1(\cdot) \asymp f_2(\cdot), $ iff
 $ f^{1/p}_1(p)/g(q) \sim f^{1/p}_2(p)/g(q).$ \par

  Note that if $ f(p) =  C \cdot h(p), \ C = const > 0, $ then
 $ f^{1/p}(p) \sim h^{1/p}(p). $ \par

 {\bf Theorem 3.1}.
 $$
 {\bf a.} \  \lambda^{-1} B_4^{1/q}(p;a,\lambda, \gamma)/g(q)
 \sim 1 + \Delta(q). \eqno(3.1)
 $$

 {\bf b.} Let us denote  $  \Lambda = \beta \max(\theta, 1/\theta) $
 for the function  $ F = F_3(p; \theta, \beta).$  We assert that:

$$
F_3(p; \theta, \beta) \asymp B_4(p; 0, \Lambda, 0); \eqno(3.2)
$$

$$
{\bf c.} \ F_3(p; \theta, \beta) \asymp G_3(p; \theta, \beta); \eqno(3.3)
$$

$$
{\bf d.}  D_4(p; a, \lambda, \gamma) \asymp B_4(p; a, \lambda, \gamma); \eqno(3.4)
$$

 {\bf Proof.}  {\sc Step 1. } We assert that $ \lambda > 0 \ \Rightarrow $

$$
\frac{\lambda^n}{n!} \le I_n(2 \sqrt{\lambda}) \le \frac{\lambda^n}{n!}
\cdot \frac{\exp(\lambda) - 1}{\lambda}. \eqno(3.5)
$$

  The first inequality is evident; let us prove the second one.
We have:
$$
\lambda^{-n} I_n(2 \sqrt{\lambda}) = \sum_{k=0}^{\infty} \frac{\lambda^k}
{k! \ (n+k+1)!} =
$$

$$
 \frac{1}{0! (n+1)!} + \frac{\lambda}{1! (n+2)!} + \frac{\lambda^2}{2! (n+3)!}
 + \ldots \le
$$

$$
\frac{1}{n!} \left( 1 + \frac{\lambda}{2!} + \frac{\lambda^2}{3!} + \ldots
\right) = \frac{1}{n!} \frac{\exp(\lambda) - 1}{\lambda}.
$$

 {\sc Step 2.}

$$
 B_4(p; a,\lambda,\gamma) \asymp \lambda^a B_3(p; \lambda, a + \gamma).
$$
 {\bf Proof.} We can assume without loss of generality that $ a $
is integer.  Further,

$$
B_4(p; a,\lambda, \gamma) \asymp \sum_{ k > a} \frac{(k-a)^p \ \lambda^k}
{\Gamma(k + \gamma + 1} = \lambda^a \sum_{n=0}^{\infty}
\frac{n^p \lambda^n}{\Gamma(k + a + \gamma + 1)}.
$$

 Let us denote for $ \alpha > 0 $ the function
$$
A(\alpha) = (\alpha + 1)^{\alpha + 0.5} \ e^{2 \alpha + 1/12}.
$$
 It follows from Stirling's formula that at $ x \ge 1 $

$$
1 \le \frac{\Gamma(x + \alpha)}{\Gamma(x) \cdot x^{\alpha}}\le
A(\alpha).
$$

 Hence

$$
A(\alpha)\ B_4(p; a,\lambda, \gamma) \asymp \sum_{m=1}^{\infty}
\frac{m^{p-a - \gamma} \ \lambda^m}{m!} = B_4(p-a-\gamma; 0,\lambda, 0).
$$
 Thus, we have obtained a simple case $ B_2(\cdot, \cdot) $
 instead of $ B_4(\cdot). $ \par

  {\sc Step 3.} We consider in this section the function $
F_3(p; \theta,\beta). $ We can suppose that $ \theta \ge 1.$
Denote

$$
F^+(p; \theta,\beta) = \sum_{k=0}^{\infty} k^p \theta^p I_k(\beta).
$$
 It follows, by virtue of the bide-side inequality $ F^+ \le F_3
\le 2 F^+ $ that $ F_3 \asymp   F^+.  $ \par
 Furthermore, we obtain using(3.5)

$$
F^+(p; \theta,\beta)  \asymp \sum_{k=1}^{\infty} \frac{k^p \ (\theta \beta)^k}
{k!} \asymp B_4(p;0,\Lambda, 0) = B_2(p; \Lambda).
$$

{\sc STEP 4.} Instead of the general case, we need to investigate
only the case $ B_2. $  We obtain using the Stirling formula with
remainder term again, $ B_2 = B_2(p; \lambda) = $

$$
\sum_{m=1}^{\infty} \frac{m^p \ \lambda^m}{m^p} \asymp \sum_{m=2}^{\infty}
\exp \left(-m \log m + m \log(e \lambda) + (p-0.5) \right).
$$

 Replacing the last sum by an integral, we can see that

$$
B_2 \asymp \int_2^{\infty} \exp( U(p;x) ) \ dx \stackrel{def}{=}
V(p) = V.
$$
$ U(p;x) = -x \log x + x \log(e \lambda)+(p-1/2) \log x. $ \par
 An integral for $ V $ can be estimated by means of classical
saddle-point method  [4], p. 119 - 127. \par
 In the next section we consider more exact estimations with non-asymptotic
 results in the particular case.\par

 \section{Probabilistic applications and non-asymptotic results.}\par

 Let $ \ p = const \ge 2, \  \{ \xi(i) \}, \ i = 1,2,\ldots, n $ be a sequence of independent centered
$ {\bf E} \xi(i) = 0 $ random variables belonging to the space $ L_p, $
i.e. such that
$$
\forall i \  \  ||\xi(i)||_p \stackrel{def}{=}  {\bf E}^{1/p}
|\xi(i)|^p < \infty. \eqno(4.0)
$$
 We denote $ \sum a(i) = \sum_{i=1}^n a(i), \ L(p) = C^p(p), $ where

$$
C(p) = \sup_{ \{\xi(i) \} } \sup_{n} \frac{||\sum \xi(i)||_p  }
{\max \left(||\sum \xi(i)||_2,
\left(\sum ||\xi(i)||^p_p \right)^{1/p} \right) }, \eqno(4.1)
$$
where $ "C" $ denotes the centered case,
 and the external $ "\sup" $ is calculated over all the sequences of
independent centered random variables that satisfy the condition
(4.0). \par
 If all the variables $ \{ \xi(i) \} $ in (4.0)  are symmetrically
distributed, independent, and  $ ||\xi(i)||_p < \infty, $ we
denote the corresponding constants (more exactly, functions of $
p) $
 $ S(p) \ ( "S" (\cdot) $ denoting the symmetrical case)
instead of $ C(p) $ and  $ \ K(p)\stackrel{def}{=} S^p(p) $
instead of $ L(p). $  Obviously, $ S(p) \le C(p), \ K(p) \le L(p).
$ It is proved in  [5]  that $ C(p) \le 2 S(p), \ L(p) \le 2^p \
K(p).$  \par
 The constant $ C(p), S(p) $ are called the exact constants in the moment
inequalities for the sums of independent random variables. They
 play very important role in the classical  probability theory
([6], 522 - 523, [7], p. 63;), in the probability theory  on the
Banach spaces ([8], [9], in the statistics and theory of Monte -
Carlo method  ([10], section 5) etc. \par
 There are many publications on the behavior of the constants $ C(p), S(p) $ as $ \ p \to \infty. $
The first estimations where obtained in [11]; Rosenthal [12]
proved, in fact, that $ C(p) \le C^p_1;  \ C_1 = const > 1; $
 here and  after $ C_j, \ j = 1,2,\ldots $ are positive finite
{\it absolute} constants.  It is proved in  [5], [13] that $ C(p)
\le 9.6 \ p / \log p. $ In the works [14], [15] are obtained the
{\it non - asymptotic } bide-sides estimations for $ S(p): $

$$
(e \sqrt{2})^{-1} \ p /\log p \le S(p) \le 7.35 \ p /\log p,
\ p \ge 2, \eqno(4.2)
$$
and there are some moment estimations for the sums independent
non-negative random variables.
 See also Latala [5], Utev [16], [17]; Pinelis and Utev [17] etc.\par
 Ibragimov R. and Sharachmedov Sh. [8], [9] and Utev [17],
 [16] have obtained the {\it explicit } formula for  $ S(p): \ S(2) = 1; $ at
$ p \in (2,4] $

$$
S(p) = \left( 1 + \sqrt{ \frac{2^p}{\pi} } \Gamma
\left( \frac{p+1}{2} \right)  \right)^{1/p};
$$

$$
p \ge 4 \ \Rightarrow S(p) = ||\tau_1 - \tau_2||_p, \eqno(4.3)
$$
 where a random variables $ \tau_j $ are independent and have the
Poisson distribution with parameters equal to 0.5: $ \ {\bf E}
\tau_j = {\bf Var} \ \tau_j = 1/2. $ \par
 As a consequence, it was obtained that as $ p \to \infty $

$$
S(p) = \frac{p}{e \cdot \log p} \left(1 + o \left(\frac{\log^2 \log p}{\log p} \right)  \right).
$$

 The following representation for the values $ L(2m), m = 2,3, 4, \ldots: $ is obtained in
 [19]:
$$
C^{2m}(2m) = L(2m) = {\bf E}(\theta - 1)^{2m} = e^{-1} \sum_{n=0}^{\infty} (n-1)^{2m} /n! = B_1(p),
$$
where the random variable $ \theta $ has the Poisson distribution
with parameter 1, and there is a hypothesis that for all values $
p \ge 4 \ C^p(p) = L(p). $ \par

  Here we denote

$$
L(p) = {\bf E}|\theta -1|^p = e^{-1} \sum_{n=0}^{\infty} |n-1|^p/n! =
  B_4(p;1,1,0), \ \ G(p) = L^{1/p}(p). \eqno(4.4)
$$
for all the values $ p \geq 4 . $ In [20]  estimations for $ C(p)
$ are obtained in the case where the sequence $ \{\xi(i) \} $ is a
sequence of martingale - differences, in  [9], [8] are presented
some generalizations for weakly dependent random variables $
\{\xi(i) \}. $ \par

 {\bf In this section we have improved the known bide-side estimations
and asymptotic for $ S(p), \ G(p) $ at $ p \to \infty, $
found the exact boundaries for the different approximations of $
S(p), \ G(p); $ describe the algorithm for the numerical
calculation of $ K(p), \ L(p); $ studied the analytical properties
of $ K(p), \ L(p) $ etc.}\par

 Note that there are many other statements of this problem: for the non-negative variables
[12], [14]; for the Hilbert space valued variables ([8], [16])
etc.
\par

   Let us denote at  $ p \ge 4: \  \delta(p) = 1/\log p, \ $

$$
  h(p) = g(p) \left(1+ \Delta(p) + \Delta^2(p) \right) =
$$

$$
 [p/(e \ \log p)] \cdot  \left(1+ \log \log p/\log p + (\log \log p/\log p)^2 \right);
$$

{\bf Theorem 4.1.}

$$
 1 = \inf_{p \ge 4} G(p)/g(p) < \sup_{p \ge 4} G(p)/g(p) = C_3,  \eqno(4.5a)
$$

{\it where}
$$
C_3 = \sup_{p \ge 4} B_2^{1/p}(p;1)/g(p) = G(C_4)/g(C_4) \approx 1.77638,
$$

$$
C_4 = \argm_{p \ge 4} B_1^{1/p}(p)/g(p) \approx 33.4610;
$$
{\it (The equality} $ C_3 \approx 1.77638 $ {\it means that } $ |C_3 - 1.77638|
\le 5 \cdot 10^{-6}); $

$$
 1 = \inf_{p =4,6,8 \ldots} C(p) /g(p) < \sup_{p=4,6,8 \ldots}
C(p)/g(p) = C_5, \eqno(4.5b)
$$

{\it where } $ C_5 = $

$$
 \inf_{p=4,6,8,\ldots} B_2^{1/p}(p;1)/g(p) = G(C_6)/g(C_6) \approx
1.77637, \ C_6 = 34;
$$

$$
1 = \inf_{p \ge 15} G(p)/h(p) < \sup_{p \ge 15}G(p)/h(p) = G(C_8)/h(C_8) =
C_{7}, \eqno(4.5c)
$$

{\it where}

$$
C_7 = \sup_{p \ge 15} B_2^{1/p}(p;1)/h(p) \approx 1.2054,
$$

$$
C_8 = \argm_{p \ge 15} B_2^{1/p}(p;1)/h(p) \approx 71.430;
$$

$$
1= \inf_{p \ge 4} S(p)/g(p) < \sup_{p \ge 4} S(p)/g(p) = C_9, \eqno(4.6a)
$$

{\it where }

$$
C_9 = \sup_{p \ge 4} W^{1/p}(p)/g(p) = S(C_{10})/g(C_{10}) \approx 1.53572,
$$

$$
 C_{10} = \argm_{p \ge 4} F_1^{1/p}(p)/g(p)  \approx 22.311;
$$

$$
1 = \inf_{p \ge 15} S(p)/h(p) < \sup_{p \ge 15} S(p)/h(p) = S(C_{12})/h(C_{12})= C_{11}, \eqno(4.6b)
$$

{\it where}

$$
C_{11}= \sup_{p \ge 15} F_1^{1/p}(p)/h(p) \approx 1.03734,
$$

$$
C_{12} = \argm_{p \ge 15} F_1^{1/p}(p)/h(p) \approx 138.149;
$$

$$
1 = \inf_{p = 16,18,20, \ldots} C(p)/h(p) < \sup_{p=16,18.20,\ldots} C(p)/h(p) =
$$

$$
C(72)/h(72) = \sup_{p = 16,18,20,\ldots} F_1^{1/p}(p;1)/h(p) \approx 1.2053. \eqno(4.6c)
$$

(We have choosen the value 15 as long as the function $ \log \log
p/\log p $ monotonically decreases at the values $ p \ge \exp(e)
\approx 15.15426 ). $ \par
  Note that our estimations and constants
 (4.5a, 4.5b, 4.5c) and (4.6a, 4.6b,4.6c) are exact and improve the constants
and estimations of Rosenthal [12]; Johnson, Schechtman, Zinn et
al. [14], [15]; Ibragimov, Sharachmedov [18],[19]; Latala  [5];
Utev [16], [17] etc. For example, $ 1/(1/\sqrt{2}) \approx
1.41421, \ 7.35 e/C_3  \approx 11.2472. $ \par

\vspace{2mm}

{\bf Theorem 4.2. }  {\it  At } $ \ p \to \infty  \ G(p) =
[p/(e \cdot \log p)] \times $

$$
 \left(1 + \frac{ \log \log p}{\log p} + \frac{1}{\log p} +
\frac{\log^2 \log p}{\log^2 p} + \frac{ \log \log p}{\log^2 p}
(1 + o(1)) \right); \eqno(4.7a)
$$

 $ S(p) = [p/(e \cdot \log p)] \times $

$$
 \left( 1 + \frac{\log \log p}{\log p} + \frac{1 - \log 2}{\log p} +
\frac{\log^2 \log p}
{\log^2 p} + o \left( \frac{\log \log p}{\log^2 p }  \right) \right).
\eqno(4.7b)
$$

\vspace{3mm}

 Let us denote by $ N = N(p), \ M = M(p) $ the
solutions of equations, which are unique:
$$
M(p) \log M(p) = p, \ \ N(p) \log(2 N(p)) = p, \eqno(4.8)
$$
for the values $ p \ge 4 $  such that $ N(p) = 0.5 M(2p). $ \par

{\bf Theorem 4.3. } {\it At } $ p \to \infty, \ m = const = 2,3,4,
\ldots,
  \ m \to \infty $
$$
G(p) = M(p)^{1-M(p)/p} \ \exp(M(p)/p) \ (1+O(\log p/p)), \eqno(4.9a)
$$

$$
C(2m) = M^{1-M(2m)/(2m)}(2m) \ \exp(M(2m)/(2m)) \ (1+O(\log m/m)),
$$

$$
S(p) = N(p) \ (e/(2N(p))^{N(p)/p} \ (1 + O(\log p/p)). \eqno(4.9b).
$$

 {\bf Theorem 4.4.} {\it Let } $ p $ {\it be even: } $ p = 2m, \ m=2,3,4,
\ldots. $ {\it Then }
$$
K(2m)=\sum_{l=0}^{2m} (-1)^l {2m \choose l} \sum_{q=0}^{2m-l}
\sum_{r=0}^l 2^{-r-q} s(2m-l,q) s(l,r),  \eqno(4.10a)
$$

$$
C^{2m}(2m) = L(2m) = \sum_{l=0}^{2m} (-1)^l {2m \choose l} \sum_{r=0}^{2m-l}
s(2m-l,r). \eqno(4.10b)
$$
{\it For integer odds values} $ p = 5,7,9, \ldots $ {\it we obtain
the representation }
$$
G^p(p) = L(p) = (2/e) + \sum_{k=0}^p (-1)^k {p \choose k} F_1(p-k). \eqno(4.11)
$$

  {\bf Proof.} First, we consider some auxiliary results. \par

{\bf 1.} {\it  In the symmetrical case, for all the values}  $ p
\in [4, \infty ) $ {\it we have: }
$$
K(p) = (2/e) \sum_{n=1}^{\infty} n^p \ I_n(1)  = F^{(S)}_1(p). \eqno(4.12)
$$
 Namely, we obtain from (4.3) for the values $ \tau_1, \ \tau_2 $
 at
$ n=1,2,\ldots: $
$$
{\bf P}(\tau_1 - \tau_2 = n) = e^{-1} \sum_{k=0}^{\infty} \frac{2^{-k} 2^{-(n+k)}}{k! \ (k+n)!} =
I_n(1)/e.
$$

{\bf 2. } On the basis of the equality (4.12)  we can offer the
numerical algorithm for $ K(p) $ investigation, calculation and
estimation. To improve convergence rate of series (4.12) we can
write:

$$
2 \pi \ I_n(1) = \int_{-\pi}^{\pi} \exp(\cos(\theta)) \ \cos(n \theta) \ d
\theta,
$$
 (see, for example, [21], p. 958, formula 5.)  After the integration
by parts, we obtain
$$
2 \pi \ I_n(1) = (-1)^m n^{-2m} \int_{-\pi}^{\pi} (\exp(\cos \theta))^{(2 m)} \
\cos (n \theta) \ d \theta,
$$
where $ m = 1,2, \ldots. $ Using the method of mathematical
induction, we conclude:

$$
(\exp \cos (\theta))^{(2m)} = \exp(\cos(\theta)) \ P_{2m}(\cos (\theta)),
$$
where $ P_{2m}(x) $ are polynomials of the degree $ 2 m $, which
can calculated by recursion

$$
P_{2m+2}(x) = (1-x^2) \left(P^{//}_{2m} + 2 P^/_{2m}(x) + P_{2m}(x)  \right) -
$$
$$
x \left(P^/_{2m}(x) + P_{2m}(x) \right)
$$
 with the initial condition $ P_0(x) = 1. $  Therefore, we obtain the following representation for
$ K(p): $

$$
\pi \ e \ K(p) = \sum_{n=1}^{\infty} n^{p-2m} \int_{-\pi}^{\pi} \exp \cos(\theta) \ P_{2m}(\theta) \
\cos (n \theta) \ d \theta.  \eqno(4.13)
$$
{\bf 3. Corollary. } {\it For even numbers } $ p = 2 m, \ m =
1,2,3, \ldots $ {\it all the numbers } $ K(p) = K(2m), L(p) =
L(2m) $ {\it are integer.} \par

 In fact, it follows from Equation (4.12) that
$$
  K(2m) =  (\pi \ e)^{-1} \ \sum_{n=1}^{\infty} \int_{-\pi}^{\pi} g^{(2m)}(\theta) \cos(n \theta)
 \ d \ \theta =
$$

$$
e^{-1} (\exp(\cos \theta))^{(2m)}(0) = (-1)^m P_{2m}(1).
$$

 It is easy to verify that all the coefficients of polynomials $ P_{2m}(x) $ are integer; thus, the
number $ P_{2m}(1) $ is integer.\par
 The second conclusion of our corollary follows from Equation (4.10), as long as all the Stirling's numbers are integer. \par

{\bf 4. For example,} $  K(6) = 31, L(6) = 41. $
 For non - integer values $ p $ we can use the method described above.
 We have obtained  using a computer program:

\begin{center}
    TABLE 1

\vspace{8mm}

\begin{tabular}{|c|c|c|c|c|c|}
\hline
p &  K(p) & L(p) & p & K(p) & L(p) \\
\hline
\hline
2 & 1 & 1 & 10.5 & 14000.4 & 41385.2  \\
\hline
4 & 4 & 4 & 11 & 30403.2 & 98253.7 \\
\hline
4.5 & 6.3358 & 6.6712 & 11.5 & 67091.3 & 236982 \\
\hline
5 & 10.4118 & 11.7358 & 12 & 150349 & 580317 \\
\hline
5.5 & 17.686 & 21.538 & 12.5 &  341951.2 & 1.44191E+006 \\
\hline
6 & 31 & 41 & 13 & 788891.0 & 3.63328E+006 \\
\hline
6.5 & 55.819 & 80.5508 & 13.5 & 1.84518E+006 & 9.27951E+006 \\
\hline
7 & 103.22 & 162.7358 & 14 & 4.37346E+006 & 2.40112E+007 \\
\hline
7.5 & 192.45 & 337.176 & 14.5 & 1.04998E+007 & 6.29176E+007 \\
\hline
8 & 379 & 715 & 15 & 2.55231E+007 & 1.66888E+008 \\
\hline
8.5 & 757.7 & 1549.28 & 15.5 & 6.27927E+007 & 4.47926E+008 \\
\hline
9 & 1126.5 & 3425.7358 & 16 & 1.56298E+008 & 1.21607 E + 009 \\
\hline
9.5 & 3015.0 & 7721.29 & 16.5 & 3.93475E+008 & 3.33839E+009 \\
\hline
10 & 6556 & 17722 & 17 & 1.00153E+009 & 9.26407E+009  \\
\hline
\end{tabular}

\vspace{4mm}

\begin{center}
  TABLE 2
\end{center}

\vspace{8mm}

\begin{tabular}{|c|c|c|}
\hline
p & K(p) & L(p) \\
\hline
\hline
17.5 & 2.57666E+009 & 2.59791E+010 \\
\hline
18   & 6.69849E+009 & 7.36008E+010 \\
\hline
18.5 & 1.75916E+010 & 2.106E + 011 \\
\hline
19 & 4.66582E+010  & 6.08476 + 011 \\
\hline
19.5 & 1.24952E+011 & 1.77473E+012 \\
\hline
20 & 3.37789E+011 & 5.22427E+012 \\
\hline
20.5 & 9.21603E+011 & 1.55177E+013 \\
\hline
21 & 2.53714E+012  & 4.64999E+013 \\
\hline
\end{tabular}

\end{center}
\vspace{8mm} {\bf 5. } Using a discrete analog of the saddle -
point method ([4], p. 262 - 264), [18]), we find that

$$
M(p) = [p/\log p] \cdot (1+ \varepsilon(p)), \
$$
where at $ p \to \infty $
$$
 \varepsilon(p) = \Delta(p) + \Delta^2(p) - \delta(p) \ \Delta(p) \ (1 + o(1)).  \eqno(4.14)
$$
Hence

$$
N(p) = [p/\log(2p)] \cdot (1+ \varepsilon(2p)) =
$$

$$
 [p/\log p] \cdot \left[1+ \Delta(p) + \Delta^2(p) - \delta(p) \ \Delta(p) \
 (1+\log 2)(1+o(1) \right].
$$

 We define the following functions and constants for the values $ p \ge P_0 = 700 $ :

$$
C_{14} = (1-\log \log P_0/\log P_0)^{-1} \approx 1.402365,
$$

$$
C_{15} = 2 \cdot \left[(1+4 \Delta^2(P_0))^{1/2} + 1 \right]^{-1} \approx 0.928958,
$$

$$
\zeta(p) = \log 2/\log(2p),
$$

$$
\varepsilon_+(p) = \Delta(p) + C_{14} \Delta^2(p), \ \varepsilon_-(p) =
\Delta(p) + C_{15} \Delta^2(p),
$$

$$
M_+ = M_+(p) = [p/\log p] \cdot (1 + \varepsilon_+(p)),
$$

$$
M_- = M_-(p) = [ p/\log p] \cdot (1+ \varepsilon_-(p)), \eqno(4.15a)
$$

$$
N_+(p) = [p/(e \cdot \log(2p))] \cdot(1+\varepsilon_+(2p)),
$$

$$
N_-(p) = [p/(e \cdot \log(2p))] \cdot(1+ \varepsilon_-(2p)). \eqno(4.15b)
$$

 More exact calculation show us that for all the values $ p \ge P_0 $
$$
 M_-(p) \le M(p) \le M_+(p), \ N_-(p) \le N(p) \le N_+(p).
$$

  Namely, one can readly see that $ \forall p \ge P_0 \ \Rightarrow $
$$
M_- \log M_- < p = M \log M < M_+ \log M_+.
$$

{\bf 6.} Let us denote

$$
 \ b_1(x,p) =  x^p/\Gamma(x+1), \ b_2(x,p) = x^p / \left(2^x \ \Gamma(x+1) \right);
$$

$$
V(x,p) =  p \log x - x \log x + x,
$$

$$
X(p) = V(M(p),p)/p = \sup_{x \ge 4} V(x,p)/p,
$$

$$
W(x,p) = p \log x - x \log x + x(1-\log 2),
$$

$$
Y(p) = W(N(p),p)/p = \sup_{x \ge 4} W(x,p)/p.
$$

 We have using the equality (4.14): $ X(p) = \log[p/(e \cdot \log p) ] + $

$$
\Delta(p) + \delta(p) + [\log(1+ \varepsilon(p)) - \varepsilon(p) + \Delta(p)
 \varepsilon(p) ] +
$$

$$
\{\delta(p) (\varepsilon(p) - \log(1+\varepsilon(p)) \} - \delta(p)
\varepsilon(p) \log(1+ \varepsilon(p))) =
$$

$$
\log[p/(e \cdot \log p)] +  X_0(p),
$$
where at  $ p \ge P_0 \ X_2(p) < X_0(p) < X_1(p), \ X_1(p)
\stackrel{def}{=} $

$$
 \Delta(p) + \delta(p) + \Delta(p) \varepsilon_+(p) + \delta(p)
[\varepsilon_+(p) -\log(1+ \varepsilon_+(p)) ], \eqno(4.16a)
$$

$$
X_2(p) \stackrel{def}{=} \Delta(p) + \delta(p) + [\log(1+\varepsilon_-(p)) - \varepsilon_-(p)] -
$$

$$
-\delta(p) \varepsilon_-(p) \log(1+\varepsilon_-(p)). \eqno(4.16b)
$$

The function $ p \to X_1(p), \ p \in [P_0, \infty) $ is monotonically decreasing and
$$
\exp(X_1(P_0)) < 1.7563, \ \lim_{p \to \infty} X_1(p) = 0. \eqno(4.16c)
$$

 In the same manner we obtain: $ Y(p) = \log[p/(e \cdot \log p)] + Y_0(p), $ where

$$
Y_0(p) = \log(1-\zeta(p)) -(1+\varepsilon(2p)) \times
$$
$$
 [1-\zeta(p) - \Delta(2p) + \delta(2p) \log(1+\varepsilon(2p)) ] + \delta(2p)
 (1+\varepsilon(2p)) \stackrel{def}{=}
$$

$$
\log g(p) + Y_0(p), \ Y_2(p) \le Y_0(p) \le Y_1(p),
$$

$$
Y_1(p) \stackrel{def}{=} \Delta(2p) + \delta(2p) +
(1+\varepsilon_+(2p)) \cdot \delta(p) \log 2/(1+\delta(p) \log 2) +
$$

$$
\varepsilon_+(2p)[\Delta(2p) + \delta(2p)],  \eqno(4.16d)
$$

$$
 Y_2(p)\stackrel{def}{=} \Delta(2p) + \delta(2p) + \varepsilon_-(2p)
 [\Delta(2p) + \delta(2p)], \eqno(4.16e)
$$

where the function $ p \to Y_1(p), \ p \in [P_1, \infty), \ P_1 =
10^6 $ is monotonically decreasing and

$$
\exp(Y_1(P_1)) < 1.442; \ Y_1(p) \downarrow 0, \ p \to \infty; \ \lim_{p \to \infty} Y_2(p) = 0.
\eqno(4.17)
$$

{\bf 7. Upper bound for } $ L(p). $ We assume in this section that
$ p \ge P_0 = 700. $ Using the well - known Stirling's formula, we
obtain for the values $ p \ge P_0, $ :

$$
 e \cdot L(p) -1.5 = \sum_{n=3}^{\infty} b_1(n-1,p) \le \int_2^{\infty} b_1(x,p) \ dx +
\sup_{x \ge 3} b_1(x,p) \le
$$

$$
(2 \pi)^{-1/2} \exp(p \cdot X(p)) + (2 \pi)^{-1/2} \int_2^{\infty} \exp (V(x,p)) \ dx.
$$
 Splitting the last integral into three parts  so that

$$
J(p) \stackrel{def}{=} \int_2^{\infty} \exp(V(x,p)) \ dx = J_1 + J_2 + J_3, \
 J_j = J_j(p), \ j = 1,2,3,
$$
 where

$$
J_1(p) = \int_2^{M - \sqrt{p}} \exp(V(x,p)) \ dx, \ J_2 = \int_{M - \sqrt{p}}^{M + \sqrt{p}}
\exp (V(x,p)) \ dx,
$$

$$
 J_3 = \int_{M + \sqrt{p}}^{\infty} \exp (V(x,p)) \ dx,
$$
we obtain for the integral $ J_2, $ taking into account the
inequalities $ M_- < M < M_+ $ and inequality: $ x \in
[M-\sqrt{p}, M+\sqrt{p} ] \ \Rightarrow $

$$
  V(x,p) \le  p X(p) - 0.5 (x-M)^2 \cdot \left( p^2 M^{-2}_+  \right) <
$$

$$
 p X(p) - 0.5(x - M)^2 \cdot p \cdot  M_+^{-2} (p):
$$
the following:
$$
J_2 \le \exp(p \cdot X(p)) \cdot \int_{M-\sqrt{p}}^{M + \sqrt{p}} \exp
\left(-0.5 p (x-M)^2 \ M^{-2}_+ \right) \ dx <
$$

$$
\exp(p \cdot X(p)) \cdot \int_{-\infty}^{\infty} \exp
\left(-0.5 p (x-M(p))^2 \  M^{-2}_+ \right) \ dx =
$$

$$
\sqrt{2 \pi} \ \exp(p \cdot X(p)) \ M_+ /\sqrt{p} \le \exp(p \cdot X(p)) \cdot \Psi_1(p),
$$
where
$$
\Psi_1(p) = \sqrt{2 \pi \ p} \ \cdot  \left[1 + \Delta +  C_{14} \Delta^2 \right]
 / \log p.
$$

 Now we estimate  the integral $ J_3. $  For the values $ x \ge M + \sqrt{p} $
the following inequalities are valid:

$$
V(x,p) \le p X(p) -   0.5 \cdot (2 p) \cdot (p/M^2_+(p)) \le p X(p) -
$$

$$
  \log^2 p \cdot (1 + \Delta + C_{14} \Delta^2 )^{-2};
$$

$$
dV(x,p)/dx \le - p /M^2_+(p) \left[x-M(p) - \sqrt{p} \right];
$$

therefore $ J_3 \le  \exp(p \cdot X(p)) \cdot \Psi_2(p), $ where $ \Psi_2(p) = $

$$
\exp \left(-\log^2 p \cdot (1 + \Delta(p) + C_{14} \Delta^2(p))^{-2} \right) \times
$$

$$
\int_{M+ \sqrt{p}}^{\infty} \exp \left(- p \ M_+^{-2}  (x-M - \sqrt{2p})   \right) \ dx  =
$$

$$
 \exp \left(-\log^2 p \cdot (1+\Delta(p) + C_{14} \Delta^2(p))^{-2} \right) \times
$$

$$
  p \cdot \left(1 + \Delta + C_{14} \Delta^2 \right)^2 \cdot \log^{-2}(p)
$$

Analogously, we find the upper estimate for $ J_1. $  \par
 Thus, $ L(p) < e^{-1} \cdot \exp(p \cdot X(p)) \times $

$$
 \left[ 1.5 \exp(-p \cdot X(p)) +
(2 \pi)^{-1/2} + \Psi_1(p) + 2 (2 \pi)^{-1/2} \ \Psi_2(p) \right] =
$$

$$
\exp(p \cdot X(p) ) \cdot \Psi^p_3(p), \eqno(4.18a)
$$
where we find by direct calculations: $ \ \Psi_3(P_0) \le 1.00826
$ and at $ p \ge P_0 $
$$
\Psi_3(p) \downarrow 1, \ p \to \infty; \ \Psi_3(p) \le 1 + C_{18} \
\log p/p.  \eqno(4.18b)
$$

{\bf 9. Lower bound for} $ L(p). $  Denote $ q = p - 1/2. $  using
 Sonin's estimate for factorials, we obtain:

$$
e L(p) \ge \sum_{n=4}^{\infty} b_1(n-1,p) = \sum_{n=3}^{\infty} b_1(n,p) \ge
$$

$$
\int_{4}^{\infty} b_1(x,p) \ dx \ge (2 \pi)^{-1/2} \ \exp(-1/12) \ \int_{4}^{\infty} \
\exp(V(x,q)) \ dx \ge
$$

$$
(2 \pi)^{-1/2} \ \exp(-1/12) \ \int_{M(q)}^{M(q) + \sqrt{ q}} \exp(V(x,q)) \ dx.
$$

 Since the following implication holds:
 $ \ x \in [M(q), M(q) + \sqrt{q} ] \ \Rightarrow $

$$
V(x,q) \ge q \ X(q)  - 0.5 (x-M(q))^2 \ q M^{-2}_-(q),
$$

we have:

$$
e L(p) \ge (2 \pi)^{-1/2} \ \exp(-1/12) \ \exp(q X(q)) \times
$$

$$
\int_{M(q)}^{M(q) + \sqrt{ q}} \exp \left[-0.5 q \ M^{-2}_-(q) \ (x-M(q))^2 \right]\ dx \ge
M_-(q) \ \times
$$

$$
0.5 \ \exp(-1/12) \ \sqrt{q} \ \exp(q X(q)) \ \left[1-\exp \left(-q^2/M^2_- \right)  \right] =
$$

$$
e \cdot \exp(p \cdot X(p)) \cdot \Psi^p_4(p),
$$
where

$$
\Psi_4(p) \downarrow 1, p \to \infty; \ \Psi_4(p) \ge 1 + C_{19} \log p/p.
$$
 Thus,

$$
\exp(p \cdot X(p)) \cdot \Psi^p_4(p) \le L(p) \le \exp(p \cdot X(p)) \cdot \Psi^p_3(p), \eqno(4.19a)
$$

$$
\Psi_3(p) \le 1 + C_{19} \log p/p, \ \Psi_4(p) \ge 1 + C_{20} \log p/p, \eqno(4.19b)
$$

$$
\Psi_{3}(p) \downarrow 1, p \to \infty; \Psi_3(P_0) \le 1.00826. \eqno(4.19c)
$$

{\bf 10. Upper and lower bounds for } $ K(p)  $ are obtained
analogously to the upper bound for $ L(p), $ but we assume in this
section that $ p \ge P_1 = 10^6. $ In brief,

$$
\sum_{k=0}^{\infty} \frac{4^{-k}}{k! \ (n+k)!} < \frac{1}{n!} \sum_{k=0}^{\infty}
\frac{4^{-k}}{k!} = \frac{\sqrt[4]{e}}{n!} < \frac{ 1.285}{n!},
$$
 hence
$$
K(p) < 2 e^{-3/4} \sum_{n=1}^{\infty} \frac{n^p \ 2^{-n}}{n!} \asymp
2 \sqrt[4]{e} \cdot B_2(p; 1/2).
$$

 Further, we conclude, using the  Stirling estimate for factorials again, that:
$$
0.5 \ e^{3/4} \ K(p) = \sum_{n=1}^{\infty} b_2(n,p) \le \int_1^{\infty} b_2(x,p) \ dx +
\sup_{x \ge 2} b_2(x,p) \le
$$

$$
 (2 \pi)^{-1/2} \exp(p \cdot Y(p)) + (2 \pi)^{-1/2} \int_{1}^{\infty} \exp (W(x,p)) \ dx.
$$
 Again, we split the last integral:

$$
 I_4 \stackrel{def}{=}  \int_2^{\infty}\exp(W(x,p)) \ dx =
$$
$$
 \left( \int_2^{N(p) - \sqrt{p}}
 + \int_{N(p)-\sqrt{p}}^{N(p) + \sqrt{p}} + \int_{N(p) +
 \sqrt{p}}^{\infty} \right) \exp \left(W(x,p) \right) dx =
$$

 $ I_5 + I_6 + I_7. $ As long as at $ x > N(p) + \sqrt{p} \ \Rightarrow $

$$
W(x,p) \le p Y(p) - 0.5p^2 N^2_+(p) = p Y(p) -
0.5 \log^2 p \cdot (1+\varepsilon_+(2p))^{-2},
$$

$$
dW/dx \le -p N^{-2}_+(x-N-\sqrt{p}),
$$
we obtain:

$$
I_7 \le \exp( p Y(p)) \cdot  p \ \log^{-2} p \
(1+\varepsilon_+(2p))^2 \times
$$

$$
 \exp \left(-0.5 \log^2 p \ (1+\varepsilon_+(2p))^{-2} \right).
$$

 Further, if $ x \in [N(p) - \sqrt{p}, N(p) + \sqrt{p}] \ $ then
$$
W(x,p) \le p Y(p) -0.5 p N^{-2}_+(p) \cdot(x-N(p))^2.
$$
 Therefore,
$$
I_6 \le \exp(p Y(p)) \cdot \sqrt{p} \ (1+\varepsilon_+(2p))/\log (2p)
$$
and $ K(p) \le 2 e^{-3/4} \ \exp(p Y(p)) \times $

$$
 \left[(2 \pi)^{-1/2}  + \sqrt{p} (1+\varepsilon_+(2p)) /\log(2p) +
2 (2 \pi)^{-1/2} p \ \log^{-2}p \right] \times
$$

$$
\left[(1+\varepsilon_+(2p))^2 \cdot \exp(-0.5 \log^2 p \
(1+\varepsilon_+(2p))^{-2} \right]. \eqno(4.20)
$$

 {\it Lower bound for } $ K(p). $ We obtain:  $ 0.5 \ e \ K(p) > $

$$
 \sum_{n=1}^{\infty} n^p 2^{-n} /n! = B_2(p;1/2) > \exp(-1/12) \
(2 \pi)^{-1/2} \ \exp(q Y(q)) \times
$$

$$
\int_{N(q)}^{N(q) + \sqrt{q} } \exp \left[ -0.5 \ q \ N^{-2}_-(q) \
(x-N(q))^2  \right] \ dx \ge
$$

$$
\exp(-1/12) \ \sqrt{\pi/2} \ q^{-1/2} \ N_-(q) \
\left(1-\exp \left(-q^2/N^2_-(q) \right) \right)/\log (2q) =
$$

$$
\exp(-1/12) \ \sqrt{\pi/2} \ \exp(q Y(q)) \ \sqrt{q} \
(1+ \varepsilon_-(2q)) \times
$$

$$
\left( 1-\exp \left(-q^2/N^2_-(q) \right) \right)/\log(2q).
$$

 Further estimations are similar to the estimation of $ L(p) $ and
can be omitted. As a result,

$$
\exp(p \cdot Y(p)) \cdot \Psi_6^p(p) \le K(p) \le \exp(p \cdot Y(p)) \cdot
\Psi^p_5(p), \eqno(4.21a)
$$
where at $ p \ge P_1 $

$$
\Psi_{5}(p) \le 1 + C_{21} \log p /p, \ \Psi_6(p) \ge 1 +
C_{22} \log p/p; \eqno(4.21b)
$$

$$
\Psi_{5}(p) \downarrow 1, \ p \to \infty;  \ \Psi_5(P_1) \le
1.000833. \eqno(4.21c)
$$

{\bf 11}. For exact computations, we need to estimate the
derivatives of our functions $ L(p), K(p). $ We show here the
estimation of  derivatives $ L^{(m)}(p), m = 1,2 \ldots. $ Namely,
$ \ e \cdot L^{(m)}(p) = $
$$
\sum_{n=3}^{\infty} \frac{(n-1)^p \ \log^m(n-1)}{n!}  \le \sum_{n=3}^{\infty}
\frac{(n-1)^p}{(n-1)!} \cdot \frac{\log^m n}{n} <
$$

$$
\sum_{n=2}^{\infty} \frac{n^p}{n!} \cdot \left( \sup_{n \ge 3}
\frac{\log^m n}{n}  \right) =
$$

$$
\left(\frac{m}{e} \right)^m \cdot \sum_{n=2}^{\infty} \frac{n^p}{n!} \asymp
\left(\frac{m}{e}\right)^m  \cdot (e B_1(p) -1). \eqno(4.22)
$$
 We estimate the derivative $ K^{(m)}(p), \ m = 1,2,\ldots $ in an
analogous way.\par
 It follows from this estimation that the functions $ L(p) $ and $ K(p) $ are infinitely
differentiable at the interval $ p \in (4, \infty). $ As long as $
L(4-0) = L(4+0) = 4, \ K(4-0) = K(4+0) = 4, $ both functions $
K(\cdot), \ L(\cdot) $ are continuous in the semiclosed interval $
[2, \infty). $ However
$$
\frac{ dK}{dp}(4-0) = \frac{dL}{dp}(4-0) \approx 3.149195,
\ \frac{dK}{dp} (4+0) \approx 3.51934,
$$
$$
\frac{dL}{dp}(4+0) \approx 3.86841,
$$
and therefore, both the functions $ K(\cdot), \ L(\cdot) $ are not
continuously differentiable in the set $ (2, \infty). $
 In the open intervals $ (2,4) $ and $ (4,\infty) $ all the
 functions $ L(p), K(p), C(p), S(p) $ are infinitely differentiable
(see (22) and [18], [19], [16] ). \par

\section{ Proof of the probabilistic results.}\par

{\bf Proof of theorem 4.1.} We find by direct calculations that $
G(C_4)/g(C_4) \approx 1.77638, $ but for we conclude  from (4.17)
that for the values $ p \ge P_0 = 700 $

$$
G(p) /g(p) \le 1.00826 \cdot 1.75913 = 1.77366,
$$
hence

$$
\argm_{p \in [4, \infty)} G(p)/g(p) \in [4, 700].
$$
 We obtain by direct calculations using known numerical methods and by
means of computer:
$$
\max_{p \in [4,700]} G(p)/g(p) = G(C_4)/g(C_4) \approx 1.77638.
$$
 Further,

$$
\inf_{p \ge 4 } G(p)/g(p) = \min \left\{\min_{p \in [4, 700]}
G(p)/g(p), \ \inf_{p > 700} G(p) /g(p) \right\}.
$$
 Our computations show that

$$
\min_{p \in [4,700]} G(p)/g(p) \approx 1.332,
$$
and it follows from (4.18a), (4.18b), (4.18c) and (4.19) that

$$
\inf_{p > 700} G(p)/g(p) = \lim_{p \to \infty} G(p)/g(p) = 1.
$$
 Thus,
$$
\inf_{p \ge 4} G(p)/g(p) = \lim_{p \to \infty} G(p)/g(p) = 1.
$$
 Analogously, $  S(C_{10})/g(C_{10}) \approx 1.53572, $ but for the
 values of
$ p \ge P_1 $ the following holds:

$$
S(p)/g(p) \le 1.0008333 \cdot 1.443 < 1.4444.
$$
 Therefore
$$
\argm_{p \in [4, \infty)} S(p)/g(p) \in [4, 1000000].
$$
We have copmuted the following:

$$
\max_{p \in [4, P_1]} S(p)/g(p) = S(C_{10})/g(C_{10}) \approx 1.53572.
$$
 Further,

$$
\inf_{p \ge 4} S(p)/g(p) = \min \{ \min_{P \in [4, P_1]}
S(p)/g(p), \inf_{p \ge P_1} S(p)/g(p) \} =
$$

$$
\lim_{p \to \infty} S(p)/g(p) = 1.
$$
 Other assertions of theorem 1 are obtained analogously. \par
{\bf Proof of theorems 4.2 and 4.3.} It follows from inequalities (4.19a),
(4.19b), (4.19c) and (4.21a), (4.21b), (4.21c) that

$$
\exp(X(p)) \cdot (1+C_{20} \log p/p) \le G(p) \le \exp(X(p)) \cdot (1+C_{19}
\log p/p), \eqno(4.23a)
$$

$$
\exp(Y(p)) \cdot (1+C_{22} \log p/p) \le S(p) \le
\exp(Y(p)) \cdot (1+ C_{21} \log p/p). \eqno(4.23b)
$$
 Substituting the expression (4.16a) and (4.16b) into equation (4.23a), we
 obtain, after simple calculation, our assertions (4.9b). We obtain (4.9a)
 in a similar way. \par

 Finally, substituting expressions (4.16d,e)  into (4.23b), we obtain
 (4.7a), (4.7b). \par

{\bf Proof of Theorem 4.4}  Since $ {\bf E} \theta_{(r)} = 1,
 r=0,1,2, \ldots, $ we conclude

$$
{\bf E} \theta^k = {\bf E} \sum_{l=0}^k s(k,l) \theta_{(r)} =  \sum_{l=0}^k  s(k,l).
$$
 Equation (4.11a) follows from the  binomial formula. The equality (4.10) is  proved analogously. \par
 Let us prove (4.11b). Since

$$
\sum_{n=2}^{\infty} n^p / n! = e \cdot B_1(p) - 1, \ p = 1,2,3, \ldots; \
\sum_{n=2}^{\infty} 1/n! = e - 2,
$$
we conclude for the values $ p = 5,7,9, \ldots:  \ e \cdot L(p) = $

$$
 1 + \sum_{n=2}^{\infty} \frac{(n-1)^p}{n!} =  1+ \sum_{n=2}^{\infty} (n!)^{-1} \cdot
\sum_{k=0}^p (-1)^k {p \choose k}  n^{p-k}  =
$$

$$
1+ \left[ \sum_{k=1}^p (-1)^k {p \choose k}
(e B_1(p-k) - 1) \right] -
$$

$$
(e B_1(0) - 2) = 2 + e \cdot \sum_{k=0}^p (-1)^k {p \choose k} B_1(p-k) -
\sum_{k=0}^p (-1)^k {p \choose k}  =
$$

$$
2+ e \sum_{k=0}^p (-1)^k {p \choose k}B_1(p-k).
$$

\section{ Some generalization and concluding remark.}

 Our results allow us to obtain some generalizations of
 symmetrically distributed random variables $ \{ \eta(i) \}, \ i =
1,2,3,\ldots.$
on the Hilbert space.  Let $ \ (H, |||\cdot|||) $ be a separable
Hilbert space with the norm $ |||\cdot|||, \ {\bf P}(\eta(i) \in
H) = 1, \ \forall i = 1,2,3 \ldots \ ||\eta||_p \stackrel{def}{=}
{\bf E}^{1/p} \left(|||\eta(i)|||^p \right) < \infty, \ p \ge 4, $

$$
Z(p) = \sup_{ \{\eta(i)\}} \sup_n \frac{||\sum \eta(i)||_p}
{\max \left( ||\sum \eta(i)||_2, (\sum ||\eta(i)||_p^p)^{1/p} \right)}.
$$
 Utev ([16], [17]) has proved that  $ Z(p) = S(p), \ p \ge 4 $ (in our notations). Therefore
$$
 1 = \inf_{p \ge 4} Z(p)/g(p) < \sup_{p \ge 4} Z(p)/g(p) = C_9 \approx 1.53572,
$$
$$
1 = \inf_{p \ge 15} Z(p)/h(p) < \sup_{p \ge 15} Z(p)/h(p)  = C_{11} \approx 1.03734.
$$
 Moreover, let us examine the following problem. For the values $ A, D \in (0, \infty),
 \ p = const \ge 4 $ we define
$$
Q^p(p, A, D) = \sup_{n \ge 1} \ \sup_{ \{ \theta(i) \} } \ {\bf E} \ |||\sum_{i=1}^n
\theta(i) |||^p,
$$
where the interior $ "\sup" \ $ is calculated over all the
sequences of $ \ H \ - $ valued independent symmetrically distributed
random variables $ \{ \theta(i) \} $ under the following conditions:

$$
\sum_{i=1}^n {\bf E} |||\theta(i) |||^2 = D^2, \ \sum_{i=1}^n
{\bf E} |||\theta(i)|||^p = A.
$$
 We denote
$$
t = 0.5 (A/D^p)^{1/(p-2)}, \
$$
then $ t  \in (0, 1/2]. $  Utev [16], [17] has also proved  that

$$
Q^p(P,A,D) = \left(A \cdot D^{-2} \right)^{1/(p-2)} \ {\bf E}|\nu_1 - \nu_2|^p,
$$
where $ \nu_1, \ \nu_2 $ are independent Poisson-distributed
random variables with parameter $ t: $

$$
{\bf P}(\nu_j = k) = \exp(-t) \cdot \ t^k /k!, \ k = 0,1,2,\ldots.
$$
Hence

$$
Q^p(p, A,D) \cdot (D^2/A)^{1/(p-2)} = {\bf E} |\nu_1 - \nu_2|^p =
$$

 $$
2 e^{-2 t} \sum_{n=1}^{\infty} n^p I_n(2 \ t) = e^{-2t} \ F_2(p;2).
 $$
 We  also denote for $ p \ge 4, \ t \in (0, 1/2] $

$$
 R(p, t) = Q(p,A,D) \cdot \left( D^2/A \right)^{1/(p(p-2))} = 2 e^{-2 t} \
\sum_{n=1}^{\infty} n^p \ I_n(2 t),
$$

$$
 u = u(t) = \sup_{p \ge 4} R(p,t)/g(p), \ T = T(t) = \argm_{p \ge 4} F(p, t),
$$
and using the above-described method, that at $ p \to \infty, \ t
= const \in (0, 1/2] $
$$
 R(p,t) = [t M(p/t)]^{1-t M(p/t)/p} \ \left(1 + o \left( \frac{\log p}{p}  \right) \right),
$$

\begin{center}
$ \ R(p, t) = [p/(e \cdot \log p) ] \times $
\end{center}

$$
\left[1 + \frac{\log \log p}{\log p} + \frac{1 + \log t}{\log p} +
\frac{\log^2 \log p}{\log^2 p} + o \left( \frac{\log \log p}{\log^2 p}
\right) \right],
$$

$$
 R^{2m}(2m,t) = \sum_{l=0}^{2m} (-1)^l {2m \choose l} \
\sum_{q=0}^{2m-l} \sum_{r=0}^l t^{q+r} \ s(2m-l,q) \ s(l,r),
$$

$$
\inf_{p \ge 4} R(p, t)/g(p) = \lim_{p \to \infty} R(p, t)/g(p) = 1;
$$
\vspace{3mm}

\begin{center}

TABLE 3 \\
\vspace{5mm}
\begin{tabular}{|c|c|c|}
\hline
 t & T(t)  & u(t)  \\
\hline
\hline
0.45 & 26.228 & 1.48566 \\
\hline
0.4  &  32.206 & 1.43438 \\
\hline
0.35 & 42.120  & 1.3815 \\
\hline
0.3 & 60.67 & 1.3281 \\
\hline
0.25 & 102.47 & 1.2732 \\
\hline
0.2 & 145.96 & 1.2163  \\
\hline

\end{tabular}

\end{center}
\vspace{2mm}

  Note that at $ t = 1/2 \ \Rightarrow T(1/2) = C_{10} \approx 22.311, \
u(1/2) = 1.53572. $ \par
 Apparently, it is interest to obtain exact constants in the
moment inequalities for sums of independent non-negative random
variables in the spirit of Refs. [12], [14], [5] etc.\par
 {\bf Aknowledgements}. We are very grateful to prof. V.Fonf, M.Lin (Ben
Gurion University, Beer - Sheva, Israel) for useful support of
these investigations. \par
 We are grateful to prof. G.Schechtman (Weizman Institute of Science, Rehovot,
Israel) for his attention.\par

\newpage

{\bf References} \\

[1] Schl\"omilch G. {\it } Zeitschrift fur Math. Und Phys., {\bf 2}, 1857, p.
140 - 154.\\

[2] Watson G.N. A Treatise on the Theory of Bessel functions.
Oxford, 1945.\\

[3] Kennet H., Rosen J. (Editor - in Chief), Michaels J.G., Gross J.L.,
at all. {\it Handbook of Discrete and Combinatorial Mathematics. } 2000,
CRC Press, Boca Raton. London, New York, Toronto, Sydney, San Francisco.\\

[4] Satshkov V.N. {\it Combinatorial Methods in discrete Mathematics.} Cambridge
University Press, 1966, Cambridge, UK.\\

[5] Latala R. {\it Estimation of Moment of Sums of independent real random
Variables.} Ann. Probab., 1997, V. 25 B.3 pp. 1502 - 1513.\\

[6] Feller W. {\it An Introduction to Probability. } V.2,  1966. John Willey
and Sons Inc., New York.\\

[7] Petrov V.V. {\it Limit Theorems of Probability Theory. Sequences of
independent Random Variables.} 1995, Oxford Science Publications. Claredon
Press, Oxford, UK.\\

[8] Talagrand M. {\it Isoperimetry and Integrability of the Sum of independent
Banach - Space valued random Variables.} Ann. Probab., 1989, V. 17 pp. 1546 - 1570.\\

[9] Dmitrovsky V.A., Ermakov S.V., Ostrovsky E.I. {\it The Central Limit
Theorem for weakly dependent Banach Space Valued Variables.} Theory Probab.
Appl., 1983, V. 28, B.1 pp. 89 - 104.\\

[10] Ostrovsky E.I. {\it Exponential Estimations for random Fields and their
Applications} (in Russian). Obninsk, OINPE, 1999.\\

[11] Dharmadhikari S., Jogdeo K. {\it Bounds on the Moment of certain
random Variables.} Ann. Math. Statist., 1969, V. 40, B.4 pp. 1506 - 1518.\\

[12] Rosenthal H.P. {\it On the subspaces of } $ L^p \ (p > 2) $ {\it spanned
by sequences of independent Variables.} Israel J. Math., 1970, V.3 pp. 273 - 253.\\

[13] Pinelis I.F. and Utev S.A. {\it Estimates of Moment of sums of
independent random variables.}  Theory Probab. Appl., 1984, V. 29 pp. 574 - 578.\\

[14] Johnson W.B., Schechtman G., Zinn J. {\it Best Constants in Moments
Inequalities for linear Combinations of independent and Changeable random
Variables.} Ann. Probab., 1985, V. 13 pp. 234 - 253.\\

[15] Johnson W.B., Maurey B., Schechtman G., Tzafriri L. {\it Symmetric
Structures in Banach Spaces.} Memoirs of the AMS., 1979, V. 217.\\

[16] Utev S.A. {\it The extremal Problems in the Moment Inequalities. } In:
Coll. Works of Siberian Branch of Academy Science USSR. Limit Theorems in
Probability Theory. 1985, V. 23 pp. 56 - 75.\\

[17] Utev S.A. {\it The extremal Problems in Probability Theory.} Probab.
Theory Appl., 1984, V. 28, B. 2 pp. 421 - 422. \\

[18] Ibragimov R., Sharachmedov Sh. {\it On the exact Constants in the Rosenthal
Inequality.} Theory Probab. Appl., 1997, V. 42 pp. 294 - 302.\\

[19] Ibragimov R., Sharachmedov Sh. {\it The exact Constant in the Rosenthal
Inequality for Sums of Independent real Random  Variables with Mean Zero.}
Theory Probab. Appl., 2001, B.1 pp. 127 - 132.\\

[20] Peshkir G., Shiryaev A.N. {\it The Khintchin Inequalities and martingale
expanding Sphere of their Actions.} Russian Math. Surveys, 1995, V.50, B.
2(305), pp. 849 - 904.\\

[21] Gradstein I.S., Ryszik I.M. {\it The Tables  of Integrals, Sums and Products.}
Academic Press, 1980, New York, London. Toronto, Sydney, San Francisco.\\

[22] Buldygin V.V., Mushtary D.M., Ostrovsky E.I., Puchalsky A.L. {\it New Trends
in Probability Theory and Statistics.} 1992, Springer Verlag, New York -
Berlin - Heidelberg - Amsterdam.\\

[23] Fedorjuk M.V. {\it The Saddle - Point Method.} Kluvner, 1990, Amsterdam,
New York.\\

[24] Figel T., Hitczenko P., Johnson W.B., Schechtman G., Zinn J.
{\it Extremal properties of Rademacher
Functions with Applications to the Khintchine and Rosenthal Inequalities.}
Trans. of the Amer.Math. Soc., V. 349 - 1027.\\

\newpage











\end{document}